\documentclass[11pt]{article}
\usepackage[letterpaper,left=1in,right=1in,top=0.9in,bottom=0.9in]{geometry}
\usepackage[T1]{fontenc}
\usepackage{lmodern}
\usepackage{amsmath,amssymb,amsthm}
\usepackage{microtype}
\usepackage{tikz}
\usepackage[section]{placeins}
\usepackage[hidelinks]{hyperref}
\theoremstyle{plain}
\newtheorem{thm}{Theorem}[section]
\newtheorem{prop}[thm]{Proposition}
\newtheorem{lem}[thm]{Lemma}
\newtheorem{cor}[thm]{Corollary}
\theoremstyle{definition}
\newtheorem{defi}[thm]{Definition}
\newtheorem{exa}[thm]{Example}
\theoremstyle{remark}
\newtheorem{rem}[thm]{Remark}
\numberwithin{equation}{section}
\title{Uniform Hyperbolicity of Automorphism\\Orbit Graphs in Rank Two}
\author{Bilal Khan\\
{\small Department of Computer Science}\\
{\small Department of Biostatistics and Health Data Science}\\
{\small Lehigh University}\\
{\small\texttt{bik221@lehigh.edu}}}
\date{}
\hypersetup{
  pdftitle={Uniform Hyperbolicity of Automorphism Orbit Graphs in Rank Two},
  pdfauthor={Bilal Khan},
  pdfsubject={Geometric group theory and automorphism orbit graphs},
  pdfkeywords={Free groups, Whitehead automorphisms, automorphic orbits,
    Schreier graphs, hyperbolicity, Bass--Serre theory}
}

\begin{document}
\maketitle

\begin{abstract}
We study automorphism orbit graphs of oriented conjugacy classes in the
free group of rank two, with every elementary Whitehead move assigned
unit cost. Collapsing signed permutations of the basis turns each orbit
component into a tree of edges and triangles. We obtain this description
by combining the known rank-two conjugacy-class stabilizer structure
with the Bass-Serre tree of the rank-two outer automorphism group. The
simplified graph and all its vertex-induced subgraphs have optimal
thin-triangle constant zero for tree components and 3/4 otherwise.
Every whole orbit component of the original graph has thin-triangle
constant at most 2, while every component of every vertex-induced
subgraph has constant at most 5/2. Both bounds are sharp in the primitive
orbit and apply to the full geometric realization, including
edge-interior triangle corners. The hereditary upper estimate also
follows from known chordality theory; the whole-component improvement
uses the geometry of adjacent permutation fibres. These results
complement earlier fixed-length and minimum-length orbit classifications.

\end{abstract}


\section{Introduction}
\label{section:introduction}

Whitehead's algorithm decides whether two words of a free group are
related by an automorphism by reducing the problem to elementary changes
of basis \cite{Whitehead1,Whitehead2}. Its analysis involves both the
length of a word representative and the cost of moving between
representatives. In the orbit graph, these become a height function,
given by cyclic word length, and a path metric, given by the number of
elementary moves. Studying the full graph reveals how the metric behaves
as representative lengths vary across an orbit.

The metric question concerns the shape of shortest move sequences: how
far can two routes between the same representatives separate, and how
much can vertex deletion enlarge those detours? Rank two is especially
accessible because its outer automorphism group has a modular-tree model
and a punctured-torus interpretation. These structures let us study whole
orbits across all cyclic lengths.

In this paper we give a global description of the orbit graphs in rank
two. Once conjugacy and signed permutations of the basis have been
factored out, each component is a tree of edges and triangles. Thus the
geometry of the entire orbit is controlled by pieces of at most three
vertices. Each fibre of the permutation quotient is a clique, so a
geodesic can visit at most two vertices of any one fibre. Combining this
constraint with the triangular block structure gives sharp uniform
bounds for the original graph as well as its permutation quotient.
The bounds are independent of the orbit and of the length of any chosen
representative, and remain uniform after arbitrary vertex deletion.

The role of stabilizers is essential. The group
$\operatorname{Out}(F_2)\cong\operatorname{GL}_2(\mathbb Z)$ acts on a
tree with finite vertex stabilizers, but quotients by subgroups can
introduce circuits of arbitrarily large length. The stabilizer of a
conjugacy class has additional structure. A filling curve has finite
stabilizer. A non-filling, nonperipheral curve has a unique disjoint
essential simple curve, and its stabilizer contains the \emph{primitive}
Dehn twist about that curve. The full twist, rather than only a proper
power, collapses the axial period responsible for large circuits. The
interaction of this surface-topological fact with the modular tree is
the mechanism behind the uniform bounds.

\subsection{Main result}

Let $F_2=\langle a,b\rangle$. The vertices of the automorphism graph
$\Omega$ are oriented conjugacy classes $[w]$ in $F_2$; two distinct
vertices are adjacent when an elementary Whitehead automorphism takes
one to the other. Let $\Pi$ be the group of eight signed permutations
of the basis, and let $\Omega^\ast$ be the quotient obtained by collapsing
each $\Pi$-orbit. Definitions and metric conventions are given in
Section~\ref{section:graphs}.

A \emph{tree of edges and triangles} is a connected graph whose blocks
are complete graphs on two or three vertices, with an isolated vertex
allowed. Equivalently, it is formed by gluing edges and triangles at
vertices so that the incidence graph of the pieces is a tree.
All edges have length one, including permutation edges in $\Omega$.
A geodesic triangle is $\delta$-thin if every point of any side lies
within distance $\delta$ of the union of the other two sides. A geodesic
metric space is $\delta$-hyperbolic here when all its geodesic triangles
are $\delta$-thin. Write $\delta(Y)$ for the infimum of such constants
for a geodesic metric space $Y$.

\begin{thm}[Uniform structure and hyperbolicity]
\label{main-theorem}
\label{tt:main}
The automorphism graphs of $F_2$ have the following properties.
\begin{enumerate}
\item Every connected component of $\Omega^\ast$ is a tree of edges
and triangles.
\item Every connected component of every vertex-induced subgraph of
$\Omega^\ast$ is $3/4$-hyperbolic in its intrinsic unit-edge metric.
More precisely, its optimal thin-triangle constant is $0$ if it is a
tree and $3/4$ otherwise.
\item Every whole connected component of $\Omega$ is $2$-hyperbolic.
Every connected component of every vertex-induced subgraph of $\Omega$
is $5/2$-hyperbolic in its intrinsic unit-edge metric. Both constants
are optimal: the primitive orbit component has constant $2$ and contains
a ten-vertex induced subgraph with constant $5/2$.
\end{enumerate}
\end{thm}

In particular, uniform hyperbolicity persists after deletion of any set
of vertices. This stability is stronger than
hyperbolicity of the original graph: deleting vertices can force new
geodesic detours. Here the separators are cliques and the quotient blocks
are edges and triangles, so those detours remain controlled. The constant
$3/4$ is optimal for $\Omega^\ast$, since triangular blocks occur; see
Example~\ref{example:triangle}. Proposition~\ref{prop:omega-sharpness}
gives explicit witnesses for the two optimal constants in $\Omega$.

\subsection{Context and proof strategy}

The study of automorphic orbits is closely connected with the complexity
of Whitehead's algorithm. Myasnikov and Shpilrain investigated orbit
sizes and the rank-two complexity of the algorithm
\cite{MyasnikovShpilrain2003}. The author's thesis \cite{Khan2003} and
the related paper \cite[Section~3.9]{Khan2004} analyse connected level neighbourhoods:
vertices of a fixed cyclic word length joined by length-preserving
elementary moves. Their principal structural result shows that
sufficiently large permutation-quotiented level neighbourhoods are chains
and gives linear bounds on their numbers of vertices. This yields a
quadratic-time analysis of Whitehead's algorithm in rank two.
Cooper and Rowland subsequently classified the finite graphs of
minimum-length representatives \cite{CooperRowland2015}. Shpilrain's
survey \cite{Shpilrain2026} reviews further algorithmic and asymptotic
developments.

The primitive orbit also has a familiar arithmetic model: its vertices
are the coprime integer pairs. Malyshev and Pak \cite{MalyshevPak2015}
consider the product replacement graph $\Gamma_2(\mathbb Z)$ on these
pairs and observe that the underlying simple graph on its positive
quadrant is a rooted binary tree.
Myropolska and Nagnibeda \cite[Section~2 and Proposition~2.4]{MyropolskaNagnibeda2014}
study the Nielsen graph $N_2(\mathbb Z)$, with additional inversion
moves, and establish its connectedness and nonamenability. In our graph
every signed permutation is a unit-cost move, and loops and parallel
edges are suppressed. These conventions are part of the exact metric
statement; the arithmetic model and its positive-quadrant tree are
already present in the earlier work.

Each whole orbit component is a simple Schreier graph of
$\operatorname{Out}(F_2)$ relative to the stabilizer of one oriented
conjugacy class. These stabilizers are elementary McCool groups.
Their general study goes back to McCool \cite{McCool1975}; see
Guirardel and Levitt \cite{GuirardelLevitt2015} for structural and
finiteness results. In particular, the rank-two finite-stabilizer and
primitive-twist alternatives used here occur in their proof of
Proposition~A.2. We give a direct surface proof because containment of
the full primitive twist determines the geometry of the quotient tree.

General methods relating graph hyperbolicity to chordality and clique
separators are well established
\cite{WuZhang2011,BermudoEdgeChordal2016,CohenEtAl2017}.
Mart\'inez-P\'erez \cite{MartinezPerez2017} relates separator and
chordality conditions to the bottleneck property and quasi-isometry to
a tree. Numerical comparisons require care: Wu and Zhang and Cohen
et al. use the four-point condition on the vertex metric, whereas the
slimness estimates of Dragan and Mohammed \cite{DraganMohammed2019}
test vertices of triangles with vertex endpoints. Our constants refer
to the Rips condition for all points of the unit-edge geometric
realization. The vertex-only bound $2$ does not by itself imply either
of our full metric bounds. However, once $6$-chordality is established,
the geometric upper bound $5/2$ follows from the edge-chordality theorem
of Bermudo et al.\ \cite[Theorem~2.4]{BermudoEdgeChordal2016}.
Remarks~\ref{remark:chordality-comparison}
and~\ref{remark:vertex-metric-comparison} give the implications precisely.

Our contribution concerns the intrinsic geometry of the full orbit
graph: every block of the permutation quotient is an edge or a triangle,
with no lower bound on word length or component size. The block theorem
gives the hereditary $5/2$ upper bound as an application of the known
chordality estimate. We sharpen the bound to $2$ for whole orbit
components and prove attainment of both constants in the primitive
orbit. To our knowledge, this global block description, the sharp
uniform whole-component bound $2$, and the realization of the hereditary
constant $5/2$ in this unit-cost Whitehead graph have not previously
been established. The $5/2$ upper estimate is thus an application of
existing graph theory, with a self-contained fibre proof included here.
The earlier level-counting results also constrain the
cyclic word-length function, information beyond the unlabelled global
block description. The value $3/4$ for a graph with only triangular
circuits is a known metric fact \cite[Theorem~3.1]{BermudoEtAl2016};
we include its elementary proof.

The proof combines punctured-torus stabilizers with the
$(2,3)$-biregular tree of the modular amalgam. A double-coset calculation
identifies orbit adjacency, and the stabilizer classification controls
the incidence quotients. Block geometry supplies the sharp simplified
bound; a clique-fibre argument treats the conjugacy-class graph and
distinguishes whole components from vertex-induced subgraphs. We recall
the classical rank-two identifications and give the stabilizer, quotient,
and metric arguments in full.

\section{The automorphism graphs and their metrics}
\label{section:graphs}

\subsection{Conjugacy classes and elementary moves}

Write $A=a^{-1}$ and $B=b^{-1}$. Elements of $F_2$ are freely reduced
words in $a,b,A,B$. A cyclically reduced word has no cancellation
between its last and first letters. The length of a conjugacy class
$[w]$ is the length of a cyclically reduced representative and is denoted
$|[w]|$. The class is oriented: conjugacy alone does not identify
$[w]$ with $[w^{-1}]$.

Every automorphism acts on conjugacy classes, since
$\phi(gwg^{-1})=\phi(g)\phi(w)\phi(g)^{-1}$. Inner automorphisms act
trivially. Let $\Pi$ be the group of signed permutations, consisting of
the maps sending $(a,b)$ to $(x,y)$ with
$x,y\in\{a,b,A,B\}$ and $y\notin\{x,x^{-1}\}$.
The elementary shifts we use are
\begin{equation}
\label{eq:basic-shifts}
\begin{aligned}
 \psi_a &: a\mapsto ab, & b&\mapsto b,\\
 \psi_b &: a\mapsto a,  & b&\mapsto ba,
\end{aligned}
\qquad
\Sigma=\{\psi_a,\psi_b,\psi_a^{-1},\psi_b^{-1}\}.
\end{equation}
Automorphisms compose as functions, with the rightmost factor acting
first. The left-hand version $a\mapsto ba$ of $\psi_a$, and the
corresponding variants of the other shifts, differ from these maps by
inner automorphisms. Consequently the eight one-letter Whitehead shifts
induce just the four maps in $\Sigma$ on conjugacy classes. Together
with the signed permutations and inner automorphisms, these are the
rank-two Whitehead generators \cite{Whitehead1,Whitehead2,LyndonSchupp}.
We also use $\Sigma$ for their images in $\operatorname{Out}(F_2)$.

\begin{defi}[Automorphism graphs]
The graph $\Omega$ has vertex set the conjugacy classes of $F_2$.
Distinct vertices $u,v$ are adjacent if $v=su$ for some
$s\in\Pi\cup\Sigma$. The graph $\Omega^\ast$ has vertex set the
$\Pi$-orbits of these classes. Two distinct orbits are adjacent when
some representatives are joined by a shift edge in $\Omega$.
For $w\in F_2$ write $w^\ast$ for its vertex in $\Omega^\ast$, and
write $\Omega_w$ and $\Omega_w^\ast$ for the corresponding connected
components.
\end{defi}

Both graphs are undirected because the generating sets are closed under
inversion on conjugacy classes. Loops and parallel edges are suppressed.
By Whitehead's generating theorem, a component is precisely an
automorphic orbit, with the indicated equivalences imposed. Every
permutation fibre of the quotient
$p:\Omega_w\to\Omega_w^\ast$ is a clique of at most eight vertices.

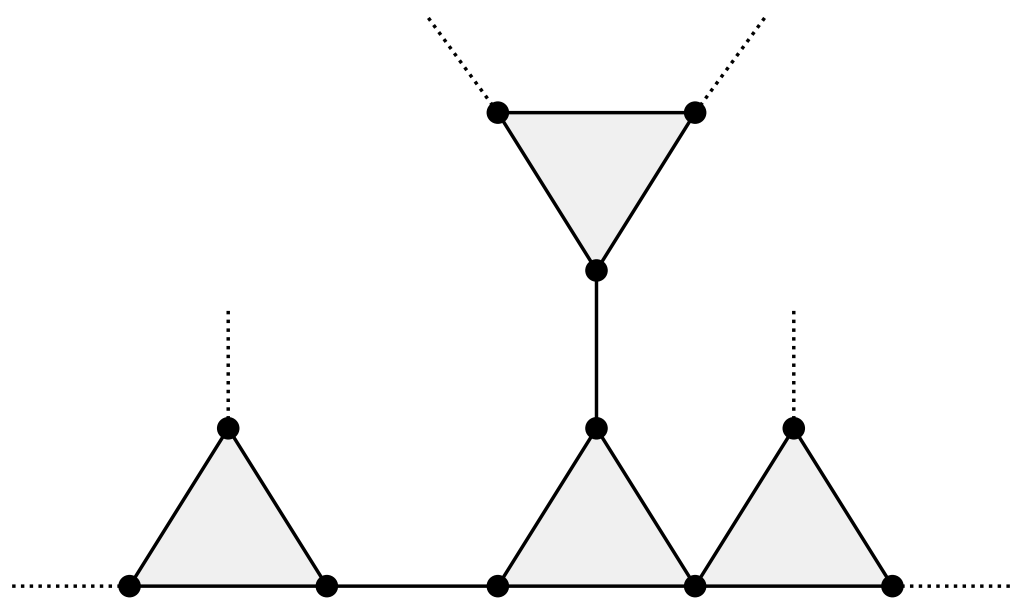
\begin{figure}[htbp]
\centering
\begin{tikzpicture}[scale=0.87,vertex/.style={circle,fill=black,inner sep=1.5pt},line width=0.65pt]
\coordinate (a) at (0,0); \coordinate (b) at (1.5,0); \coordinate (c) at (.75,1.2);
\coordinate (d) at (-1.3,0); \coordinate (e) at (-2.05,1.2); \coordinate (f) at (-2.8,0);
\coordinate (g) at (2.25,1.2); \coordinate (h) at (3,0);
\coordinate (i) at (.75,2.4); \coordinate (j) at (0,3.6); \coordinate (k) at (1.5,3.6);
\fill[gray!12] (a)--(b)--(c)--cycle;
\fill[gray!12] (d)--(e)--(f)--cycle;
\fill[gray!12] (b)--(g)--(h)--cycle;
\fill[gray!12] (i)--(j)--(k)--cycle;
\draw (a)--(b)--(c)--cycle (a)--(d) (d)--(e)--(f)--cycle;
\draw (b)--(g)--(h)--cycle (c)--(i) (i)--(j)--(k)--cycle;
\draw[densely dotted] (f)--(-3.7,0) (e)--(-2.05,2.1)
 (h)--(3.9,0) (g)--(2.25,2.1) (j)--(-.55,4.35) (k)--(2.05,4.35);
\foreach \v in {a,b,c,d,e,f,g,h,i,j,k} \node[vertex] at (\v) {};
\end{tikzpicture}
\caption{A schematic tree of edges and triangles. The cliques meet at
cut vertices, and their incidence graph is a tree; dotted segments
indicate possible continuations.}
\label{fig:block-structure}
\end{figure}

\subsection{Path metrics and vertex deletion}

Every edge of $\Omega$ and $\Omega^\ast$ has length one. We regard each
connected component as a geodesic metric graph, including the interior
points of edges. Thus the bounds in Theorem~\ref{main-theorem} concern
thin triangles for arbitrary points, not only graph vertices.
For a vertex set $V_0$, the induced subgraph retains exactly the edges
with both endpoints in $V_0$; each connected component receives its own
intrinsic metric. In particular, deleting a vertex ball means deleting
the vertices in that ball and all their incident edges.

\begin{prop}[Comparison of the two metrics]
\label{prop:metric-comparison}
For vertices $u,v$ in one component of $\Omega$,
\[
 d_{\Omega^\ast}(p(u),p(v))\leq d_\Omega(u,v)
 \leq d_{\Omega^\ast}(p(u),p(v))+1.
\]
\end{prop}
\begin{proof}
Projection collapses permutation edges and sends a shift edge to an
edge or a vertex, giving the first inequality. On conjugacy classes,
conjugation by $\Pi$ permutes $\Sigma$; this follows by interchanging
the basis letters and inverting either one. Therefore every shift edge
between two permutation classes can be lifted starting at any chosen
representative of its initial class. Explicitly, if
$s\pi_1u=\pi_2v$ with $s\in\Sigma$ and $\pi_1,\pi_2\in\Pi$, then
\[
 (\pi_1^{-1}s\pi_1)u=\pi_1^{-1}\pi_2v,
 \qquad \pi_1^{-1}s\pi_1\in\Sigma
\]
as maps of conjugacy classes. Lift a shortest path in $\Omega^\ast$
one edge at a time, starting at $u$. The lifted path has the same
length and ends in the fibre of $v$. At most one permutation edge
then joins it to $v$.
\end{proof}

On the vertex set of $\Omega$, giving permutation edges weight zero
produces a pseudometric whose metric quotient is the vertex set of
$\Omega^\ast$ with its path metric. At the level of geometric
realizations, the quotient can retain parallel edges; the definition of
$\Omega^\ast$ suppresses these before assigning unit lengths. We keep
the unit-edge metric throughout, so that the finite permutation fibres
remain visible in $\Omega$. All graph-theoretic terms below refer to
the simple graphs just defined.

\section{Conjugacy-class stabilizers}
\label{tt:section-stabilizers}

The alternatives below specialize the known rank-two McCool-group
structure \cite[Proposition~A.2 and its proof]{GuirardelLevitt2015}
to one oriented conjugacy class. We record the surface argument and
the precise twist-normalizer inclusion needed in Section~\ref{tt:section-tree-quotients}.

We use the classical rank-two identification
\[
 G=\operatorname{Out}(F_2)\cong\operatorname{GL}_2(\mathbb Z).
\]
It is induced by abelianization: its kernel on $\operatorname{Aut}(F_2)$
is $\operatorname{Inn}(F_2)$. It also identifies $G$ with the extended
mapping class group of the once-punctured torus $S$. We use homeomorphisms of the punctured surface; no boundary component
is fixed pointwise. In particular there is no boundary-twist kernel. These standard identifications, simple
curve slopes, and Dehn twists are discussed in
\cite{LyndonSchupp,FarbMargalit}.
Concretely, every integral unimodular matrix acts linearly on the torus
and fixes the puncture at the origin. Its outer action has the prescribed
abelianization, so the rank-two identification ensures that this realizes
the action of every element of $G$ on free homotopy classes.

An essential simple curve on $S$ is one that is neither null-homotopic nor
peripheral. Its unoriented isotopy class corresponds to a primitive vector
in $\mathbb Z^2$ modulo multiplication by $-1$. Distinct such classes
intersect and fill $S$. A closed curve is \emph{filling} if its geometric
intersection number with every essential simple curve is positive.
The trivial class and the powers of the peripheral class are not filling.
Here and below intersection number includes the multiplicity of a power.

\begin{lem}[Bounded intersections with a filling curve]
\label{tt:bounded-intersections}
If $c$ is a filling closed curve on $S$, then for every $n\geq 0$ the set
\[
 \mathcal S_n(c)=\{\delta:\delta\text{ is an unoriented essential simple
 curve on }S,\ i(c,\delta)\leq n\}
\]
is finite.
\end{lem}
\begin{proof}
Write $c=u^m$ with $m\geq1$ and $u$ not a proper power. It suffices to
prove the assertion for $u$, since $i(c,\delta)=m\,i(u,\delta)$.
Choose a complete finite-area hyperbolic metric and take the geodesic
representative $C$ of this root. Its image is a finite embedded graph,
allowing vertices at its self-intersections. Since $c$ fills, the
complementary regions are disks and once-punctured disks. One may either
choose a generic metric or subdivide at all multiple intersections.

Put a simple curve $\delta$ in minimal position with $C$ and avoid the
vertices of $C$. Geodesic representatives realize geometric intersection
numbers; a small perturbation at multiple intersections gives this normal
position without changing the number counted with multiplicity. Cutting
$\delta$ at its intersections with $C$ gives at most $n$ embedded arcs
in the complementary regions.

There are finitely many choices of the edges containing the at most $2n$
endpoints and their orders on each edge. For any fixed choice there are
finitely many isotopy classes of disjoint arc systems in each face. In a
disk these are determined by their noncrossing pairing. In a
once-punctured disk there are still only finitely many: after specifying
the pairing, one also specifies which complementary region contains the
puncture. Equivalently, an embedded arc with both endpoints on the outer
boundary has only the two possibilities for the side containing the
puncture; arbitrary multiple winding is incompatible with embeddedness.
The arcs have no endpoints at the puncture.

Gluing these finitely many possibilities across the edges of $C$ gives
finitely many isotopy classes of $\delta$. Since $c$ fills, no essential
simple curve is disjoint from $C$, so there is no unrecorded closed
component lying in a face. This proves the assertion.
\end{proof}

\begin{lem}[Finite stabilizer in the filling case]
\label{tt:filling-stabilizer}
If $[w]$ is filling, then
$H_w=\operatorname{Stab}_G([w])$ is finite.
\end{lem}
\begin{proof}
Choose $n$ large enough that $\mathcal S_n(w)$ contains two distinct
slopes. The group $H_w$ permutes this finite set because it preserves
intersection number with $w$. An element of the kernel fixes the two
unoriented slopes. If $v_1,v_2$ are primitive vectors representing them,
its matrix must send $v_j$ to $\epsilon_jv_j$, where
$\epsilon_j\in\{1,-1\}$. The two vectors are linearly independent over
$\mathbb Q$, so the two signs determine the matrix. The kernel thus has
at most four elements, and the image is a finite permutation group.
\end{proof}

\begin{prop}[The stabilizer alternatives]
\label{tt:stabilizers}
For an oriented conjugacy class $[w]$ exactly one of the following
geometric cases applies.
\begin{enumerate}
\item If $w=1$, then $H_w=G$ and the orbit is a single vertex.
\item If $[w]$ is peripheral and nontrivial, then $H_w$ is the
determinant-one subgroup $\operatorname{SL}_2(\mathbb Z)$, and its
component in $\Omega^\ast$ is a single vertex.
\item If $[w]$ is filling, then $H_w$ is finite.
\item Otherwise there is a unique unoriented essential simple curve
$\delta$ with $i(w,\delta)=0$. If $t_\delta$ is a primitive Dehn twist
about it, then
\[
 \langle t_\delta\rangle\ \leq H_w\ \leq
 N_\delta:=\operatorname{Stab}_G(\delta),
 \qquad
 N_\delta/\langle t_\delta\rangle\cong C_2\times C_2.
\]
In particular, the entire primitive twist subgroup, not merely a
finite-index subgroup of it, is contained in $H_w$.
\end{enumerate}
\end{prop}
\begin{proof}
The first assertion is immediate. A peripheral class is $[c^m]$, where
$c=aba^{-1}b^{-1}$ and $m\ne0$. A mapping class sends it to
$[c^{m\det(g)}]$. A nontrivial element of a free group is not conjugate
to its inverse: an element conjugating it to its inverse would reverse
its axis in the free Cayley tree, hence fix a vertex or invert an edge,
which is impossible for the free action. Thus the stabilizer is exactly
the determinant-one subgroup. A signed permutation of determinant $-1$
interchanges the two possible conjugacy classes, so their permutation
quotient is one vertex. The filling case is
Lemma~\ref{tt:filling-stabilizer}.

In the remaining case some essential simple curve $\delta$ is disjoint
from $w$. If two distinct such curves existed, take their geodesic
representatives together with that of $w$. The two simple curves
intersect and fill the once-punctured torus; their complement consists
of disks and a punctured disk. A closed curve disjoint from both must
therefore be trivial or peripheral, contrary to the hypotheses.
This also handles powers of simple curves: distinct essential simple
curves on $S$ have positive intersection. Hence $\delta$ is unique.

Every element fixing $[w]$ preserves the unique class $\delta$.
Conversely, a representative of $w$ can be chosen outside an annular
support for $t_\delta$, so this twist fixes its oriented free homotopy
class. If $w$ is a power of $\delta$, take a parallel representative.
This proves the two inclusions.

Conjugate so that $\delta$ has slope $(1,0)$. In the column convention
for abelianization, its unoriented stabilizer and a primitive twist are
\begin{equation}
\label{tt:cusp}
 N=\left\{
 \begin{pmatrix}\epsilon&n\\0&\eta\end{pmatrix}:
 \epsilon,\eta\in\{1,-1\},\ n\in\mathbb Z\right\},
 \qquad
 R=\begin{pmatrix}1&1\\0&1\end{pmatrix}.
\end{equation}
The diagonal-sign map $N\to C_2\times C_2$ has kernel
$\langle R\rangle$. Moreover $gRg^{-1}=R^{\epsilon\eta}$ for
$g\in N$. This proves the last assertion as well as normality of the
primitive twist subgroup.
\end{proof}

\section{A modular-tree model}
\label{tt:section-quotient}

Use the following matrices in $G$:
\begin{equation}
\label{tt:matrices}
 Q=\begin{pmatrix}0&-1\\1&0\end{pmatrix},\qquad
 U=\begin{pmatrix}0&-1\\1&1\end{pmatrix},\qquad
 S_0=\begin{pmatrix}0&1\\1&0\end{pmatrix}.
\end{equation}
Put
\[
 K=\langle Q,S_0\rangle=\Pi,\qquad
 D=\langle U,S_0\rangle,\qquad
 C=\langle -I,S_0\rangle.
\]
Their orders are $8,12,4$, respectively, and $K\cap D=C$.
The standard modular-group presentation gives
\begin{equation}
\label{tt:amalgam}
 G=K\mathop{*}_C D.
\end{equation}
For an explicit check on the embeddings, the presentation is
\[
 \langle Q,U,S_0\mid Q^4=1,\ Q^2=U^3,\ S_0^2=1,
       \ S_0QS_0=Q^{-1},\ S_0US_0=U^{-1}\rangle.
\]
Indeed $\operatorname{SL}_2(\mathbb Z)$ has presentation
$\langle Q,U\mid Q^4=1,Q^2=U^3\rangle$, and adjoining $S_0$
gives the determinant-$-1$ coset with the displayed conjugation action.
This is also the amalgam presentation of the indicated dihedral groups.
See \cite{Serre}; the rank-two amalgam is explicitly recalled in
\cite[p.~3]{Zimmermann2021}.

Let $T$ be its Bass--Serre tree. Its vertex sets are $G/K$ and $G/D$,
and an edge $gC$ joins $gK$ to $gD$. Vertices of the two types have
degrees two and three. The action preserves the types and has no
edge inversions. We retain the two types throughout the proof.
The central element $-I$ acts trivially on $T$; faithfulness is not
required in any of the quotient arguments below.

\begin{defi}[The incidence construction]
\label{tt:incidence-definition}
For a bipartite multigraph $Z$ with vertex types $K$ and $D$, let
$\mathcal I(Z)$ be the simple graph on its $K$-vertices in which two
distinct vertices are adjacent precisely when they share a $D$-neighbour.
Parallel incidence edges in $Z$ have no effect on $\mathcal I(Z)$.
\end{defi}

Write $X=\mathcal I(T)$. Thus $X$ is the tree of triangles obtained by
replacing each $D$-star by a clique, or the line graph of the trivalent
tree obtained by suppressing the $K$-vertices.

\begin{lem}[The double-coset adjacency]
\label{tt:double-coset-adjacency}
For distinct $gK,g'K\in G/K$,
\[
 gK\sim_Xg'K\quad\Longleftrightarrow\quad
 g^{-1}g'\in KUK=KRK,
\]
where $R$ is the matrix in~\eqref{tt:cusp}, representing $\psi_b$.
\end{lem}
\begin{proof}
The three neighbours of $D$ in $T$ are $K,UK,U^2K$.
Here $U^3=-I\in C$. The element $S_0\in C$ exchanges the last two
neighbours, since $S_0US_0=U^{-1}$. Consequently
$D\setminus C=CUC$. A $D$-vertex adjacent to $gK$ has the form
$gkD$, $k\in K$, and its other $K$-neighbours are $gkdK$ with
$d\in D\setminus C$. This proves the first equality of relations.
Direct multiplication gives $R=Q^{-1}U$, so $KRK=KUK$.
\end{proof}

\begin{lem}[Quotients commute with the incidence construction]
\label{tt:quotient-incidence}
For every subgroup $H\leq G$, there is a natural simple-graph isomorphism
\[
 H\backslash X\cong\mathcal I(H\backslash T).
\]
In $H\backslash T$ the vertices and edges are orbits; parallel incidence
edges are retained until applying $\mathcal I$.
\end{lem}
\begin{proof}
Both vertex sets are the $H$-orbits of the $K$-vertices of $T$.
An adjacent pair upstairs certainly gives a shared $D$-neighbour orbit
downstairs. Conversely, suppose two $K$-vertex orbits have incident edges
to the same $D$-vertex orbit. Choose lifts of the two incidence edges,
with $D$-endpoints $d_1,d_2$. Some $h\in H$ sends $d_2$ to $d_1$.
Translate the second edge by $h$. The two lifted $K$-endpoints now
share the actual vertex $d_1$ and hence are adjacent in $X$. Their
orbits are unchanged. Suppressing coincident endpoints gives exactly
the stated simple-graph isomorphism.
\end{proof}

\begin{prop}[Identification of the orbit graph]
\label{tt:orbit-identification}
For $H=H_w$, the component $\Omega_w^\ast$ containing $w^\ast$ satisfies
\[
 \Omega_w^\ast\cong H\backslash X
       \cong\mathcal I(H\backslash T).
\]
\end{prop}
\begin{proof}
The orbit of $[w]$ is the left $G$-set $G/H$ via $gH\mapsto g[w]$.
Permutation classes are therefore $K\backslash G/H$. Inversion gives
the bijection
\begin{equation}
\label{tt:inversion}
 KgH\longmapsto Hg^{-1}K
\end{equation}
with the vertex set $H\backslash G/K$ of $H\backslash X$.

Let $\Sigma=\{\psi_a^{\pm1},\psi_b^{\pm1}\}\subset G$.
The two positive shift matrices are
\[
 L=\begin{pmatrix}1&0\\1&1\end{pmatrix}=S_0RS_0,
 \qquad R=\begin{pmatrix}1&1\\0&1\end{pmatrix}.
\]
Conjugation by signed permutations interchanges these matrices and
their inverses, so $K\Sigma K=KRK=KUK$.

Two permutation classes $Kg_1H,Kg_2H$ are joined by a shift edge
exactly when $g_2\in K\Sigma K g_1H$. Taking inverses, and using
$\Sigma^{-1}=\Sigma$, gives
$g_2^{-1}\in H g_1^{-1}K\Sigma K$. By
Lemma~\ref{tt:double-coset-adjacency}, this is precisely adjacency of
the corresponding $H$-orbits in $X$. The second isomorphism follows
from Lemma~\ref{tt:quotient-incidence}.
\end{proof}

\section{The structure of orbit quotients}
\label{tt:section-tree-quotients}

\begin{lem}[Rooted tree quotients]
\label{tt:rooted-quotient}
If a group acts on a tree fixing a vertex $o$, its orbit multigraph is
a tree. In particular the quotient by a finite group acting without
inversions is a tree.
\end{lem}
\begin{proof}
Distance from $o$ is constant on vertex orbits. Every non-root vertex
has a unique predecessor one level nearer $o$, and this predecessor
is equivariant. A non-root vertex orbit consequently has exactly one
edge orbit to the preceding level. This also rules out parallel edges
to that level. Edges only join consecutive levels, so connectedness
and the unique-predecessor property prove the assertion.

For a finite group, the convex hull of a finite orbit is a finite
invariant subtree. Its centre is a vertex or an edge. In the latter
case absence of inversions implies that the endpoints are fixed.
Thus a fixed vertex exists.
\end{proof}

\begin{lem}[A primitive twist in the Bass--Serre tree]
\label{tt:primitive-twist}
Suppose $\langle R\rangle\leq H\leq N$, with $N$ as in~\eqref{tt:cusp}.
After identifying parallel edges, $H\backslash T$ is a tree.
\end{lem}
\begin{proof}
Since $R=Q^{-1}U$, we have $R^{-1}K=U^{-1}K$. Thus $R^nD$
is adjacent to both $R^{n-1}K$ and $R^nK$, for every $n\in\mathbb Z$.
The vertices $R^nK$ are pairwise distinct: otherwise a nonzero power
of the unipotent matrix $R$ would belong to the finite group $K$.
The same argument with the finite group $D$ shows that the vertices
$R^nD$ are pairwise distinct. These vertices and incidences therefore
form an embedded bi-infinite line $\mathcal A$ in the tree $T$.
The element $R$ translates this line by two edges, so it is hyperbolic
with axis $\mathcal A$ and translation length two.
Hence $\langle R\rangle\backslash\mathcal A$
consists of one $K$-vertex and one $D$-vertex joined by two parallel
edges.

Every $K$-vertex on the axis has degree two, so has no off-axis
branches. Each $D$-vertex on the axis has exactly one off-axis branch.
The translates of any one such branch are all the off-axis branches
along $\mathcal A$, and distinct branches are disjoint away from their
attachment points. Thus $Z_0=\langle R\rangle\backslash T$ is the
two-edge circuit just described, with one rooted tree attached at
its $D$-vertex.

The subgroup $\langle R\rangle$ is normal in $N$, and $N/\langle
R\rangle$ has order four. Each element of $N$ preserves the unique
axis of $R$, since it conjugates $R$ to $R$ or $R^{-1}$. Therefore
$H/\langle R\rangle$ acts on $Z_0$, fixes both vertices of its axial
circuit (the types are different), and preserves the unique attached
rooted tree. The quotient of that rooted tree is again a tree by
Lemma~\ref{tt:rooted-quotient}. The two axial edges are either retained
or identified with one another. After suppressing their multiplicity,
the resulting graph is a tree. Taking the two successive quotients is
the same as taking $H\backslash T$.
\end{proof}

\begin{rem}[The primitive twist condition]
\label{tt:full-twist-essential}
The distinction between a primitive twist and a proper power has a
geometric effect. For $H=\langle R^m\rangle$, $m\geq4$,
the axis quotient has $m$ distinct $K$-vertices and $m$ distinct
$D$-vertices. The incidence construction contains an induced circuit
of length $m$ on the axial $K$-vertices. Proposition~\ref{tt:stabilizers}
excludes precisely this obstruction for nonperipheral conjugacy-class
stabilizers: whenever a disjoint curve exists, the primitive twist
itself fixes the class.
\end{rem}

\begin{lem}[Incidence graphs of bipartite trees]
\label{tt:incidence-tree}
Let $Z$ be a bipartite tree with at least one $K$-vertex, whose
$K$-vertices have degree at most two
and whose $D$-vertices have degree at most three. Then $\mathcal I(Z)$
is a tree of edges and triangles, allowing a singleton.
\end{lem}
\begin{proof}
Each $D$-vertex with at least two neighbours contributes the clique
on those neighbours, of size two or three. Two distinct such cliques
share at most one vertex; otherwise $Z$ would have a circuit of length
four. Delete the $D$-leaves from $Z$. The remaining incidence graph of
these cliques and the $K$-vertices is a tree (or a single vertex).
In particular no cyclic chain of distinct cliques is possible. Each
of the contributed cliques is a block, and these are all the blocks.
This proves the assertion.
\end{proof}

\begin{proof}[Proof of the structural assertion in Theorem~\ref{tt:main}]
For the trivial and peripheral classes the graph is a singleton by
Proposition~\ref{tt:stabilizers}. For a filling class, $H_w$ is finite,
so $H_w\backslash T$ is a tree by Lemma~\ref{tt:rooted-quotient}.
For a non-filling, nonperipheral class, conjugating the unique disjoint
curve to slope $(1,0)$ conjugates $H_w$ to a subgroup between
$\langle R\rangle$ and $N$. Conjugation of subgroups gives isomorphic
quotient trees: $H\backslash T\to gHg^{-1}\backslash T$ is induced
by $v\mapsto gv$. Lemma~\ref{tt:primitive-twist} therefore shows that
the simple underlying graph of $H_w\backslash T$ is a tree.

At any quotient vertex, incidence edge orbits can be represented by
edges at one chosen lift, so its degree, counted with multiplicity,
does not exceed the degree of that lift. Thus the simple quotient has
$K$-degrees at most two and $D$-degrees at most three. Parallel
incidences do not change $\mathcal I$. Apply
Proposition~\ref{tt:orbit-identification} and
Lemma~\ref{tt:incidence-tree} to conclude that $\Omega_w^\ast$ is
a tree of edges and triangles.
\end{proof}

\section{Uniform hyperbolicity and vertex deletion}
\label{tt:section-metric}

The exact constant for the simplified graph follows directly from its
block structure. The next lemma is the relevant special case of the
small-constant characterization in \cite[Theorem~3.1]{BermudoEtAl2016}.

\begin{lem}[The sharp constant for triangular block graphs]
\label{lemma:sharp-blocks}
Let $Y$ be a tree of edges and triangles with unit-length edges.
Every connected component $Y_0$ of every vertex-induced subgraph of $Y$
is again a tree of edges and triangles, and
\[
 \delta(Y_0)=
 \begin{cases}
  0,&Y_0\text{ is a tree},\\
  3/4,&Y_0\text{ contains a triangle}.
 \end{cases}
\]
\end{lem}
\begin{proof}
Deleting vertices restricts each block to a smaller clique. These
cliques still meet only at cut vertices and their incidence graph is a
forest. Each resulting connected component therefore has the stated
block structure. It suffices to compute the constant for $Y$ itself.

Every component outside a block $B$ attaches to $B$ at a single vertex;
otherwise the block incidence graph would contain a circuit. Thus $B$
is isometrically embedded, and a geodesic cannot leave $B$ and then
return to it. The geometric diameter of a triangular block is $3/2$,
and that of an edge block is $1$.

Let $\alpha$ be a side from $x$ to $y$ of a geodesic triangle, and let
$P$ be the union of its other two sides. For $p\in\alpha$, choose a
block $B$ containing $p$ and a nontrivial portion of $\alpha$; a constant
side needs no argument. The intersection $\alpha\cap B$ is a geodesic
segment with endpoints $r,s$ and length at most $3/2$. Each endpoint is
either $x$ or $y$, or the unique attachment vertex through which the
corresponding endpoint $x$ or $y$ is connected to $B$. Consequently every
$x$--$y$ path, including $P$, contains both $r$ and $s$. It follows that
\[
 d(p,P)\leq\min\{d(p,r),d(p,s)\}
       \leq\tfrac12\operatorname{length}(\alpha\cap B)\leq\tfrac34.
\]
This proves the upper bound for arbitrary points of the metric graph.

If $Y$ is a tree, its triangles are tripods and $\delta(Y)=0$.
Otherwise it contains an isometrically embedded triangular block,
which is a circle of length three. Choose antipodal points $x,y$ on
this circle, take one semicircle as $\alpha$, and choose a third point
on the other semicircle for the triangle's third vertex. The midpoint
of $\alpha$ has distance $3/4$ from the other two sides. The embedding
is isometric, so this distance is also $3/4$ in $Y$. This proves the
lower bound.
\end{proof}

For $\Omega$ the clique geometry of the fibres gives sharper information
than their cardinalities. The following lemma allows arbitrary edges
between adjacent fibres; it therefore continues to apply after deletion.

\begin{lem}[Clique fibres over triangular block graphs]
\label{lemma:clique-fibres}
Let $Z$ be a tree of edges and triangles. Suppose a connected unit-edge
graph $Y$ has a map $f:V(Y)\to V(Z)$ whose nonempty fibres are cliques,
and every edge of $Y$ has endpoints in one fibre or in adjacent fibres.
Then $\delta(Y)\leq5/2$, including for triangles with edge-interior
endpoints. Moreover, $\delta(Y)\leq2$ if
\begin{equation}
\label{eq:adjacent-fibre-distance}
 d_Y(u,v)\leq2\quad\text{whenever }f(u)\sim_Z f(v).
\end{equation}
\end{lem}
\begin{proof}
We first record the attachment property used below. Let
$S\subseteq V(Z)$ be nonempty and connected, let $U$ be a component
of $Z\setminus S$, and put
\[
 \partial_S U=\{b\in U:b\text{ has a neighbour in }S\}.
\]
All edges between $U$ and $S$ lie in a single block $B$. To see this,
join the respective endpoints of any two distinct such edges by simple
paths inside $S$ and $U$. Their union with the two edges is a simple
circuit, hence lies in one block of $Z$. Since every edge belongs to
just one block, the assertion follows. Thus $\partial_S U=V(B)\setminus S$
has at most two vertices. Every simple path in $U$ between boundary
vertices stays in $B$: an alternative path between two distinct such
vertices would, together with their edge in $B$, give a circuit outside
that block. In particular the path uses at most two vertices. If $U$
has two distinct neighbours in $S$, then $B$ is their triangle and
$\partial_S U$ consists of its third vertex.

Let $\alpha$ be a nonconstant side from $x$ to $y$ of a geodesic
triangle and $P$ the union of its other two sides. Constant sides
need no argument. Choose a simple $x$--$y$ path $Q\subseteq P$;
such a path exists because $P$ is a finite union of edge segments.
If $Q$ contains no graph vertex, it has length at
most one, so $\alpha$ also has length at most one. The same conclusion
holds if $\alpha$ contains no graph vertex. In either case its points
are within $1/2$ of $\{x,y\}\subseteq Q$. We may therefore assume
both paths contain vertices. Let $S$ be the set of images under $f$
of the vertices on $Q$. It is nonempty and connected in $Z$. Every
vertex in a fibre indexed by $S$ has distance at most one from $Q$.

A geodesic visits at most two vertices of any fibre, and these occur
consecutively: the subpath between two distinct vertices in a fibre
has length at most one. Consequently, projecting the vertex sequence
of $\alpha$ and suppressing consecutive repetitions gives a simple
path $\gamma$ in $Z$.

Consider a maximal run of $\gamma$ outside $S$ with both a predecessor
and a successor in $S$. These two vertices of $S$ are distinct, since
$\gamma$ is simple. The attaching-block observation shows that the
run consists of one vertex, the third vertex of a triangle. Let $r$
be the last vertex of $\alpha$ in its preceding fibre and $s$ the
first vertex in its succeeding fibre. The subpath from $r$ to $s$
has length $L\leq3$:
it uses two edges between fibres and at most one edge inside the
intermediate fibre. For every point $p$ on this subpath,
\begin{equation}
\label{eq:interior-run}
 d_Y(p,Q)\leq\min\{d_Y(p,r)+1,d_Y(p,s)+1\}
             \leq L/2+1\leq5/2.
\end{equation}
Here $d_Y(p,r)+d_Y(p,s)=L$, since $\alpha$ is geodesic.
Under~\eqref{eq:adjacent-fibre-distance}, the preceding and succeeding
fibres are adjacent in $Z$, so geodesicity gives $L\leq2$ and the
same estimate is at most $2$.

Suppose next that an initial run of $\gamma$ outside $S$ is followed
by a vertex of $S$. Let $u$ be the first graph vertex of $\alpha$.
Since $f(u)\notin S$ and $x\in Q$, we have $u\ne x$, so $x$ is
interior to an edge. The path $Q$ must follow that edge towards
its other endpoint $v$; otherwise $u$ would be on $Q$. Hence $v\in Q$
and $f(u)\in\partial_S U$ for the outside component $U$ containing
the run. The last fibre index of the run also lies in $\partial_S U$,
so the attachment
property allows at most two fibres. Write $t=d_Y(x,u)\in(0,1)$.
If there were two outside fibres, their
attaching block would have just one vertex in $S$, namely $f(v)$.
The subpath to its first return vertex $r$ in that fibre would have
length at least $t+2$. The alternative through $v$ has length at most
$(1-t)+1=2-t$, a contradiction. Thus there is only one outside fibre.
The subpath from $x$ to its first vertex $r$ in a fibre indexed by $S$
has length $L\leq t+2<3$. Since $x\in Q$ and $d_Y(r,Q)\leq1$,
every point of it has distance at most $(L+1)/2<2$ from $Q$.
Indeed, at distance $a$ along this subpath from $x$, the bound is
$\min\{a,L-a+1\}\leq(L+1)/2$.
Terminal runs are treated in the same way.

If all vertices of $\alpha$ project outside $S$, they lie in one
component $U$ of $Z\setminus S$. As in the preceding argument, both
$x$ and $y$ are interior to edges whose other endpoints lie on $Q$.
The first and last vertices of $\gamma$ belong to $\partial_S U$.
The attachment property therefore shows that $\gamma$ uses at most
two fibres, which are adjacent if distinct. Between its first and
last vertices, $\alpha$ has length at most three. Adding its two
partial end-edges gives length at most five, so every point is within
$5/2$ of $\{x,y\}\subseteq Q$. If
\eqref{eq:adjacent-fibre-distance} holds, those first and last vertices
have distance at most two when they belong to different fibres, and
at most one otherwise. Geodesicity then bounds the total length by
four and the distance to $\{x,y\}$ by two.

All remaining full edges of $\alpha$ have both endpoint fibres indexed
by $S$. A nearest endpoint puts each of their points within $3/2$ of
$Q$. A remaining partial end-edge lies within distance one of $x$ or $y$.
Together these cases show that $\alpha$ lies within $5/2$ of $Q$,
or within $2$ under~\eqref{eq:adjacent-fibre-distance}. Since
$Q\subseteq P$, the lemma follows.
\end{proof}

\begin{rem}[Comparison with general chordality bounds]
\label{remark:chordality-comparison}
The graphs in Lemma~\ref{lemma:clique-fibres} are $6$-chordal: they
have no induced circuit of length greater than six. Indeed, an induced
circuit of length at least four uses at most two vertices of a fibre,
consecutively. If it uses at least three fibres, their cyclic sequence
is a simple circuit in $Z$, hence a triangle. A circuit using at most
two fibres has at most four vertices. This proves the assertion.

In the terminology of \cite{BermudoEdgeChordal2016}, a graph is
$(k,m)$-edge-chordal when every circuit of length at least $k$ has a
shortcut consisting of one edge of length at most $m$. The parameters
$k,m$ are nonnegative real numbers. Since circuit lengths are integers
in a unit-edge graph, $6$-chordality implies
$(6+\varepsilon,1)$-edge-chordality for every $\varepsilon>0$.
Their Theorem~2.4, stated for the full geometric Rips constant, gives
\[
 \delta(Y)\leq 1+\frac{6+\varepsilon}{4}
       \quad\text{for every }\varepsilon>0,
 \qquad\text{hence}\qquad \delta(Y)\leq\frac52.
\]
This applies directly to our locally finite orbit graphs and their
induced subgraphs. Thus the hereditary $5/2$ upper bound is already a
consequence of the block structure and this known theorem.
Lemma~\ref{lemma:clique-fibres} supplies a self-contained proof and uses
the additional adjacent-fibre distance condition to give the sharper
bound $2$ for whole components. Proposition~\ref{prop:omega-sharpness}
establishes attainment of both constants in the orbit graphs themselves.

For comparison, Dragan and Mohammed
\cite[Corollary~4]{DraganMohammed2019} give the vertex-slimness bound
$\lfloor6/4\rfloor+1=2$ for finite $6$-chordal graphs. A finite induced
subgraph containing a given triangle with graph-vertex corners preserves
its chosen geodesic sides, so this vertex estimate applies in our infinite
graphs as well. It is a different consequence of the same chordality
hypothesis.
\end{rem}

\begin{rem}[Vertex slimness and edge interiors]
\label{remark:vertex-metric-comparison}
Let $s(G)$ test only graph vertices on geodesic triangles whose three
corners are graph vertices. If $s(G)\leq2$, then every point on a side
of such a triangle lies within $1/2$ of a side vertex, and therefore
within $5/2$ of the other sides. This argument retains the restriction
on triangle corners.

Allowing arbitrary corners requires more information. For example,
$s(C_{11})=2$: each vertex-to-vertex geodesic has at most five edges,
so every vertex on it lies within two of one endpoint. Equality occurs
on the triangle whose successive side lengths are $5,3,3$. On the
geometric circle of circumference eleven, every geodesic has length
at most $11/2$, giving $\delta(C_{11})\leq11/4$. Conversely, take
antipodal points $x,y$ and a third corner on one of the two semicircles.
The midpoint of the other semicircle has distance $11/4$ from the
remaining sides. Thus
\[
 s(C_{11})=2,\qquad \delta(C_{11})=11/4>5/2.
\]
Every proper induced subgraph of $C_{11}$ is a forest, so even a
hereditary vertex-slimness bound $2$ does not imply the full metric
bound $5/2$. The $6$-chordality used in the preceding remark is a
stronger hypothesis. Moreover, the induced graph in
Proposition~\ref{prop:omega-sharpness} is $6$-chordal with vertex
slimness $2$ and geometric constant $5/2$; hence $6$-chordality alone
cannot supply the whole-component bound $2$.
\end{rem}

\begin{proof}[Proof of the metric assertions in Theorem~\ref{tt:main}]
The block description has been proved in Section~\ref{tt:section-tree-quotients}.
Lemma~\ref{lemma:sharp-blocks} gives the bound $3/4$ for $\Omega^\ast$
and the asserted exact constants, including after arbitrary vertex
deletion.

Let $p:\Omega_w\to\Omega_w^\ast$ be the permutation quotient.
Every fibre $F_v=p^{-1}(v)$ is a clique. A permutation edge remains
inside a fibre, and a shift edge projects to an edge or a vertex.
Proposition~\ref{prop:metric-comparison} gives
\eqref{eq:adjacent-fibre-distance} for the whole component. Apply
Lemma~\ref{lemma:clique-fibres} to obtain the bound $2$.

Now let $Y$ be a connected component of any vertex-induced subgraph
of $\Omega_w$. The nonempty intersections $F_v\cap V(Y)$ are still
cliques, and edges still project to edges or vertices of the same
tree of edges and triangles $\Omega_w^\ast$. The first assertion of
Lemma~\ref{lemma:clique-fibres} gives $\delta(Y)\leq5/2$. No lifting
of quotient edges through deleted vertices is used here. The lower
bounds, and hence optimality of both constants, are proved in
Proposition~\ref{prop:omega-sharpness}.
\end{proof}

\section{Geodesics, examples, and the role of rank two}
\label{section:examples}

\begin{cor}[Shortest paths and induced circuits]
\label{cor:geodesics}
Between two vertices in a component of $\Omega^\ast$ there is a unique
shortest vertex path. Every simple circuit in a component has length
three. These conclusions hold in every connected component of every
vertex-induced subgraph of $\Omega^\ast$.
\end{cor}
\begin{proof}
The block incidence tree determines the sequence of blocks traversed
between the two vertices. Within each block the entry and exit vertices
are joined by their unique edge. This gives the unique shortest vertex
path. A simple circuit cannot traverse a cut vertex twice, so it lies
in a single block and is a triangle. Deleting vertices restricts each
clique to a smaller clique and preserves the tree arrangement.
\end{proof}

\begin{cor}[Ball complements]
For every $r\geq0$ and every vertex $u$, each connected component
obtained by deleting the vertex ball of radius $r$ about $u$ is
$3/4$-hyperbolic in $\Omega^\ast$ and $5/2$-hyperbolic in $\Omega$,
respectively, with its intrinsic path metric.
\end{cor}
\begin{proof}
Both complements are vertex-induced subgraphs, so
Theorem~\ref{main-theorem} applies.
\end{proof}

\begin{exa}[An elementary triangle]
\label{example:triangle}
Let $Q$ be the signed permutation $a\mapsto b$, $b\mapsto A$.
In $\operatorname{Out}(F_2)$ one has
\begin{equation}
\label{eq:triangle-relation}
 \psi_a\psi_b^{-1}\psi_a=Q.
\end{equation}
Indeed, with the representatives in~\eqref{eq:basic-shifts}, the
composite sends $a$ to $abA$ and $b$ to $A$. It is the composition
of $Q$ with the inner automorphism $g\mapsto aga^{-1}$.
The induced maps on conjugacy classes therefore satisfy the same
identity without an inner factor. For instance,
\[
 [abb]\xrightarrow{\psi_a}[abbb]
 \xrightarrow{\psi_b^{-1}}[bAbAb]
 \xrightarrow{\psi_a}[AAb]=Q[abb].
\]
The three permutation classes have positive representatives
$[abb]$, $[abbb]$, and $[ababa]$ of respective lengths $3,4,5$.
They form a genuine geodesic triangle after the permutation quotient.
Thus triangles persist, even though all larger simple circuits are
excluded. Their bounded size is part of the structure, rather than
an effect of restricting to short words.
\end{exa}

\begin{exa}[Peripheral and primitive classes]
The classes of $1$ and of $(aba^{-1}b^{-1})^m$, $m\ne0$, give singleton
components in $\Omega^\ast$. For $[a^m]$, $m\ne0$, the stabilizer has
the concrete form
\[
 H_{a^m}=\left\{
 \begin{pmatrix}1&n\\0&\eta\end{pmatrix}:
 n\in\mathbb Z,\ \eta\in\{1,-1\}\right\}.
\]
To see this, fixing $[a^m]$ forces its nonzero homology vector
$(m,0)$ to be fixed. Conversely, the displayed matrices are generated
by $R$ and the signed permutation $a\mapsto a$, $b\mapsto B$,
both of which fix $[a^m]$. This is an explicit instance of the
primitive-twist alternative in Proposition~\ref{tt:stabilizers}.
\end{exa}

\begin{prop}[Optimal constants in the primitive orbit]
\label{prop:omega-sharpness}
The orbit component $\Omega_a$ of the primitive class $[a]$ satisfies
$\delta(\Omega_a)=2$. It contains a ten-vertex induced subgraph $H$
with $\delta(H)=5/2$.
\end{prop}
\begin{proof}
We first identify the vertices of $\Omega_a$ with the primitive integer
vectors $(p,q)$, where $\gcd(p,q)=1$, using abelianization. Every such
vector is the first column of an element of
$\operatorname{GL}_2(\mathbb Z)$. The stabilizer of $(1,0)$ consists
of the matrices
\[
 \begin{pmatrix}1&n\\0&\eta\end{pmatrix},
 \qquad n\in\mathbb Z,\quad\eta\in\{1,-1\}.
\]
They are generated by $R$ and $\operatorname{diag}(1,-1)$, both of
which fix $[a]$. Conversely, fixing $[a]$ fixes its homology vector.
Thus the two stabilizers agree, proving that abelianization is a
bijection on this orbit. In these coordinates, permutation adjacency
means that the unordered pairs of absolute coordinates agree, while
the four shift images are
\begin{equation}
\label{eq:primitive-moves}
 (p,q)\longmapsto(p+q,q),\ (p-q,q),\ (p,p+q),\ (p,q-p).
\end{equation}

Take the ten vertices
\begin{equation}
\label{eq:extremal-vertices}
\begin{array}{c|ccccc}
 i&0&1&2&3&4\\ \hline
 v_i&(-2,-1)&(-3,-1)&(-4,-1)&(-1,-4)&(3,-4)\\[2pt]
 i&5&6&7&8&9\\ \hline
 v_i&(3,-7)&(-3,7)&(-3,4)&(-3,1)&(-2,1)
\end{array}
\end{equation}
The adjacency formulas give exactly the ten-cycle
$v_0v_1\cdots v_9v_0$ and the two additional edges $v_1v_8,v_4v_7$
in their induced subgraph $H$. Explicitly, its permutation edges have
index pairs
\[
 \{0,9\},\ \{1,8\},\ \{2,3\},\ \{4,7\},\ \{5,6\},
\]
and its shift edges have index pairs
\[
 \{0,1\},\ \{1,2\},\ \{3,4\},\ \{4,5\},\
 \{6,7\},\ \{7,8\},\ \{8,9\}.
\]
There are no other edges, by~\eqref{eq:primitive-moves}. The labels
\[
 5,4,3,2,1,0,1,2,3,4.
\]
on $v_0,\ldots,v_9$ differ by at most one on every edge and are
realized by paths to $v_5$. They are therefore the exact distances
from $v_5$ in $H$.
Consequently the paths
\[
 \alpha=v_0v_1v_2v_3v_4v_5,\qquad
 \beta=v_5v_6v_7v_8v_9,\qquad
 \gamma=v_9v_0
\]
are geodesic sides of a triangle in $H$. Put $P=\beta\cup\gamma$.
The distance labels from $V(P)$, in the same vertex order, are
\[
 0,1,2,2,1,0,0,0,0,0:
\]
again they are $1$-Lipschitz on every edge and realized by paths.
Since $P$ is a union of full edges, a path from a vertex outside $P$
first enters $P$ at a vertex. Thus both $v_2$ and $v_3$ have metric
distance two from $P$, including its edge interiors. Every path from the midpoint
of $v_2v_3$ to this union first passes through $v_2$ or $v_3$, so that
midpoint has distance $5/2$ from the other sides. The upper bound
already proved gives $\delta(H)=5/2$.

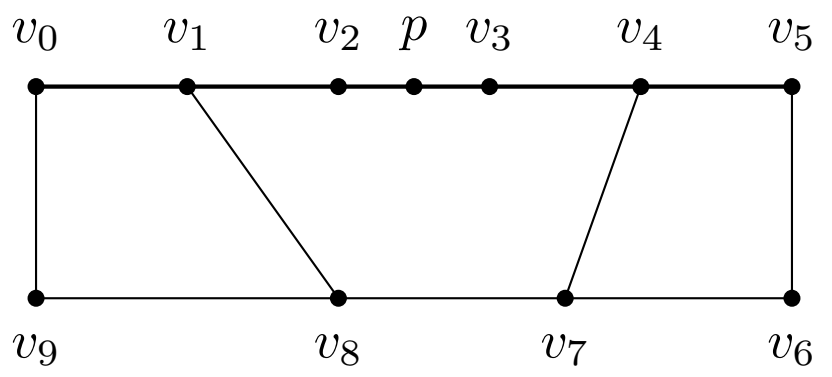
\begin{figure}[htbp]
\centering
\begin{tikzpicture}[x=1cm,y=0.7cm,every node/.style={font=\small}]
 \coordinate (v0) at (0,1); \coordinate (v1) at (1,1);
 \coordinate (v2) at (2,1); \coordinate (v3) at (3,1);
 \coordinate (v4) at (4,1); \coordinate (v5) at (5,1);
 \coordinate (v6) at (5,-1); \coordinate (v7) at (3.5,-1);
 \coordinate (v8) at (2,-1); \coordinate (v9) at (0,-1);
 \draw[thick] (v0)--(v1)--(v2)--(v3)--(v4)--(v5);
 \draw (v5)--(v6)--(v7)--(v8)--(v9)--(v0);
 \draw (v1)--(v8) (v4)--(v7);
 \foreach \i in {0,...,5} {
   \fill (v\i) circle (1.6pt);
   \node[above=3pt] at (v\i) {$v_{\i}$};
 }
 \foreach \i in {6,...,9} {
   \fill (v\i) circle (1.6pt);
   \node[below=3pt] at (v\i) {$v_{\i}$};
 }
 \fill (2.5,1) circle (1.6pt);
 \node[above=3pt] at (2.5,1) {$p$};
\end{tikzpicture}
\caption{The ten-vertex induced subgraph $H$, with all edges of length
one. The top path is a geodesic of length five from $v_0$ to $v_5$;
the other two sides run through $v_9$. Its midpoint $p$ has distance
$5/2$ from their union.}
\label{fig:omega-sharp}
\end{figure}

For the whole component, also use the primitive vectors
\[
 r=(-4,3),\qquad s=(-7,3).
\]
The following paths in $\Omega_a$ have lengths $4,4,1$:
\begin{equation}
\label{eq:whole-extremal-triangle}
 v_0v_1v_2rs,\qquad sv_6v_7v_8v_9,\qquad v_9v_0.
\end{equation}
To check geodesicity in the whole graph, index each permutation fibre
by the increasing pair of absolute coordinates. The quotient path
\[
 (1,2)\;--\;(1,3)\;--\;(3,4)\;--\;(3,7)
\]
is induced, by~\eqref{eq:primitive-moves}. In a tree of edges and
triangles an induced path is the unique quotient geodesic. Any path
of length three from either $v_0$ or $v_9$ to $s$ would therefore have
to lift precisely these three quotient edges. The shift formulas give
the only such lifts from the respective starting vertices:
\[
 \begin{aligned}
 (-2,-1)&\longmapsto(-3,-1)\longmapsto(-3,-4)
                   \longmapsto(-3,-7),\\
 (-2,1)&\longmapsto(-3,1)\longmapsto(-3,4)
                   \longmapsto(-3,7).
 \end{aligned}
\]
Neither ends at $s=(-7,3)$. Thus both length-four paths
in~\eqref{eq:whole-extremal-triangle} are geodesics, and the third
side is an edge.

The vertex $v_2=(-4,-1)$ has no neighbour on the other two sides of
this triangle. Their vertices have different unordered absolute
coordinate pairs, and the four shift images of $v_2$ are
\[
 (-5,-1),\quad(-3,-1),\quad(-4,-5),\quad(-4,3),
\]
none of which lies on those sides. Their union consists of full edges,
so a path from $v_2$ first enters it at a vertex. Since $v_2$ is not
on them either, its distance from their union is at least two, and the path
$v_2v_1v_0$ realizes distance two. Hence $\delta(\Omega_a)\geq2$.
The upper bound for whole components proves equality.
\end{proof}

The argument uses two features specific to rank two. The modular group
has an amalgam with finite factors whose incidence degrees are two and
three, and a non-filling, nonperipheral curve on the punctured torus has
only one disjoint essential simple-curve class. The first feature bounds
the local pieces; the second controls the possible periods introduced
by a stabilizer. Together they explain why the same finite block geometry
occurs throughout every automorphic orbit, despite the variation in
cyclic word lengths and in the sizes of minimum-length levels.

\section*{Acknowledgements}
The author thanks Alexei Myasnikov for helpful discussions, comments,
and suggestions. Research support from the Mathematics Department of
the City University of New York Graduate Center, the Center for
Computational Sciences at the Naval Research Laboratory, and the
Advanced Engineering and Sciences division of ITT Industries is
gratefully acknowledged.

\bibliographystyle{plain}
\bibliography{references}
\end{document}